# Modulated MPC for Arm Inductor-less MVDC MMC with Reduced Computational Burden

Sandro P. Martin, *Member, IEEE*, Hui Li, *Fellow, IEEE*,
and Olugbenga M. Anubi, *Member, IEEE*

***Abstract*--A modulated model predictive controller is designed for an inductor-less modular multilevel converter targeting an MVDC solid-state transformer application. The underlying optimization problem is formulated such that a unique closed-form solution is derived from a highly accurate dynamic system model, which significantly reduces computation time. Unlike other closed-form methods, the proposed controller achieves inductor-less current control, has a reduced sampling speed enabled by a novel model-based circulating energy compensator, and does not require tuned PI controllers to generate current references. Simulation and experimental results demonstrate that the proposed controller achieves transient power flow control and steady state circulating energy control despite the low system inertia.**

## I. Introduction

THE medium-voltage DC (MVDC) application field is expanding, and includes electric ships, datacenters, electric railways, and the integration of distributed generation and renewable energy resources to the grid [1-4]. The MVDC solid-state transformer (SST), a power electronics-enabled technology that fills the equivalent role of the traditional AC transformer, enables these applications. One popular converter for implementing the MVDC SST is the modular multilevel converter (MMC), which gained its popularity in HVDC applications. By connecting two MMCs in a back-to-back or front-to-front configuration, they may perform DC-AC-DC or AC-DC-AC conversion. The MMC offers several important advantageous features such as modularity, low harmonics, obviation of DC bus capacitors, and good fault withstanding capability [5-7].

Despite these beneficial characteristics, the MMC suffers from large-sized cell capacitors and arm inductors. Recently, a number of research efforts have been reported to reduce the size of these passive components [8-10]. The small cell capacitor, arm inductor-less MMC can achieve high power density, but its control becomes more challenging due to stability issues inherent to a low-inertia system, strong coupling between AC and DC stages, and increased circulating energy, resulting in higher loss and requiring faster capacitor voltage balancing control.

A classical PI regulator-based control system is usually adopted for MMCs with large cell capacitors and arm inductors [11-14]. This approach may use *dq*-based PI controllers for the power flow control, and PI or P controllers for the circulating energy control and capacitor voltage balancing control. However, this traditional control method has several issues when it is applied to the MMC with small passive components. First, the lack of a linear time-invariant (LTI) mathematical model that accurately describes the coupling effect for this MIMO system prevents *a priori* design, instead requiring extensive parametric tuning. Although methods such as [14-15] propose techniques for designing such controllers *a priori*, they are typically of limited usefulness: [14] designs a passivity-based controller from a model subjected to a multi-frequency orthogonal coordinate transformation, but its derivation is highly complicated; [15] designs a controller from a small-signal model based on harmonic state-space analysis, but is also complicated and based on an MMC with large passive components. In addition to the parametric tuning requirement, the traditional PI-based controller typically suffers from low bandwidth, which cannot achieve the required stability and circulating energy suppression of a low inertia system [8].

Within the last decade, model predictive control (MPC) has gained popularity in MMC applications. MPC does not require a time-invariant system model, is conceptually simple, and has inherent MIMO capabilities. In MMC applications, two popular MPC strategies are the optimal switching state (OSS) MPC and the optimal voltage level (OVL) MPC, also known as modulated MPC. The OSS MPC technique determines the switching state of the MMC directly, resulting in a very fast dynamic control. However, the calculation burden of this process is exponentially related to the number of cells in the MMC, meaning that it requires a very fast controller. By contrast, the modulated MPC method is more comparable to the PI approach, as its output is a voltage reference, which can then be used for modulation. The computational burden of the typical modulated MPC approach is only linearly related to the number of MMC cells, but its response time is slower than the OSS MPC. An overview of existing OSS and OVL MPC techniques has been presented in [16].

In addition to their separate disadvantages, both the OSS and OVL techniques tend to suffer from two shared drawbacks: a

The authors are with the Center for Advanced Power Systems at the Florida State University, Tallahassee, FL 32310 USA
E-mails: spmartin2@fsu.edu, li@eng.famu.fsu.edu, oanubi@fsu.edu

variable switching frequency, and the need to optimize weighting factors in the cost function manually. However, a variation on the modulated technique first described in 2018 that is based partly on the deadbeat control principle may eliminate both of these disadvantages, and may furthermore remove the dependency of the computational burden on the number of cells present in the MMC. This technique was demonstrated in [17] for an MMC based on phase-shifted PWM, and [18] for an MMC based on level-shifted PWM. In these cases, the deadbeat-like modulated MPC formulation is such that it solves both the power flow control and the circulating energy control in one step. The algorithm has a unique solution, and does not require weighting factor selection because the exact control requirements are already included in the predicted optimal references. Thus, this variation on modulated MPC can be implemented with very low computational burden and has the additional benefit of being insensitive to the number of cells in the MMC.

The fast pace of MPC research for MMC has resulted in numerous improvements over the state of the art in the year since the publication of [16]. A region-less explicit MPC (RL-EMPC), which solves the OVL MPC optimization problem offline while a fast controller solves a simpler online computation without exponential memory requirements, was demonstrated in [21] with an FPGA controller implementation. In [22], an OSS direct model predictive power control (DMPPC) technique is applied to perform the power, circulating current, and voltage balancing controls without requiring PI controllers for reference generation. In [23], an OVL-based MPC is proposed to reduce the computational cost, where the control output is solved by Diophantine equations rather than through minimization of a cost function. A machine learning controller trained on an OSS MPC is presented in [24] to achieve a similar control quality as the MPC with a reduced computational burden. Another variation on the OVL MPC is proposed in [25] based on a submodule grouping technique to reduce the computational burden.

However, even in MPC-based MMCs, passive component sizes tend to remain large [17-25]. Table A-1 in the Appendix shows a per-unit comparison of the sizes of the passive components in a sample of both PI-based and MPC-based MMCs reported so far, as well as their sampling frequencies. The proposed MPC in this paper is shown in the last row for comparison, where it is readily apparent that its inductance is the smallest among the MMCs surveyed, and its capacitance is among the smallest. While other surveyed MMCs have smaller per-unit capacitances, they offset their small capacitor sizes with larger inductors.

The existing literature on controller design for SSTs with small passive components is limited. The authors of [8] have introduced the DC-DC modular soft-switching solid-state transformer (M-S4T). The M-S4T is a modular converter that uses small passive components and experiences tight coupling. The authors of [8] compare the performance of a PI-based controller and a deadbeat-type controller in balancing the submodule voltages, balancing the magnetizing currents, maintaining the load voltage, and tracking the power reference. The results indicate that the PI-based controller performs significantly worse than the deadbeat controller, and fails to achieve full output power. This suggests that a deadbeat or other model-predictive controller is more appropriate for controlling a low-inertia, tightly coupled converter.

In this paper, a modulated MPC is designed for an arm inductor-less MMC targeted for a MVDC transformer application. The controller in this paper differs from existing methods in several ways, leading to the following novel contributions:

- MPC based on a dynamic model that explicitly defines the circulating energy;
- MPC achievement of inductor-less current control;
- MPC operation with a reduced sampling speed via a novel model-based circulating energy compensator;
- MPC with model-based external PI controller to generate inductor-less current references.

Although some of the above have been achieved in other MPCs, this research is the first to 1) combine all into a single simple controller and 2) do so for a low-inertia system. The rest of this paper is organized as follows. The MMC dynamic model is introduced in Section II. The modulated MPC and novel feed-forward circulating energy compensator are developed in Section III. In Section IV, a simulated MVDC application case study is provided to validate the controller and compare it against the controller in [17]. In Section V, experimental results on a downscaled prototype are presented to verify the proposed control. A conclusion is presented in Section VI.

## II. Mathematic Model of the Arm Inductor-less MMC with Phase-Shifted Square Wave Modulation

The topology of the MMC-based MVDC SST is shown in Fig. 1(a). Due to its symmetric feature, only one side of the DC-DC converter modeling and control is examined in this paper, as the analysis can extend to the whole converter. The one side of this DC-DC converter, which is an MMC, is shown in Fig. 1(b). By adopting a phase-shifted square wave modulation (PS-SWM), the MMC can achieve an arm inductor-less feature, so $L_{arm}$ and $L_{dc}$ of Fig. 1 (b) are parasitic arm inductances and a DC cable inductance, respectively. The modulation method was described in detail in [10], and Fig. 1(c) shows the modulation waveforms. To generate PWM signals, phase-shifted triangular carriers with period $T_{sw}$ are compared against a square wave reference $m = (d_{dc} - 1) \pm d_{ac}$, where $d_{dc}$ and $d_{ac}$ are DC and AC control variables defined in (1),

$$\begin{cases} d_{dc} = \dfrac{V_{dc}}{v_c^*} \\ d_{ac} = \left|\dfrac{v_{ac}^*}{v_c^*}\right| \times square(\omega t) \end{cases} \tag{1}$$

and where $square(\omega t)$ is a square wave with range (-1, +1) and period $T = 2\pi/\omega$. The DC reference, $v_c^*$, is the arm capacitor voltage magnitude; the AC reference, $v_{ac}^*$, is the magnitude of the AC voltage; and $V_{dc}$ is the DC input voltage.

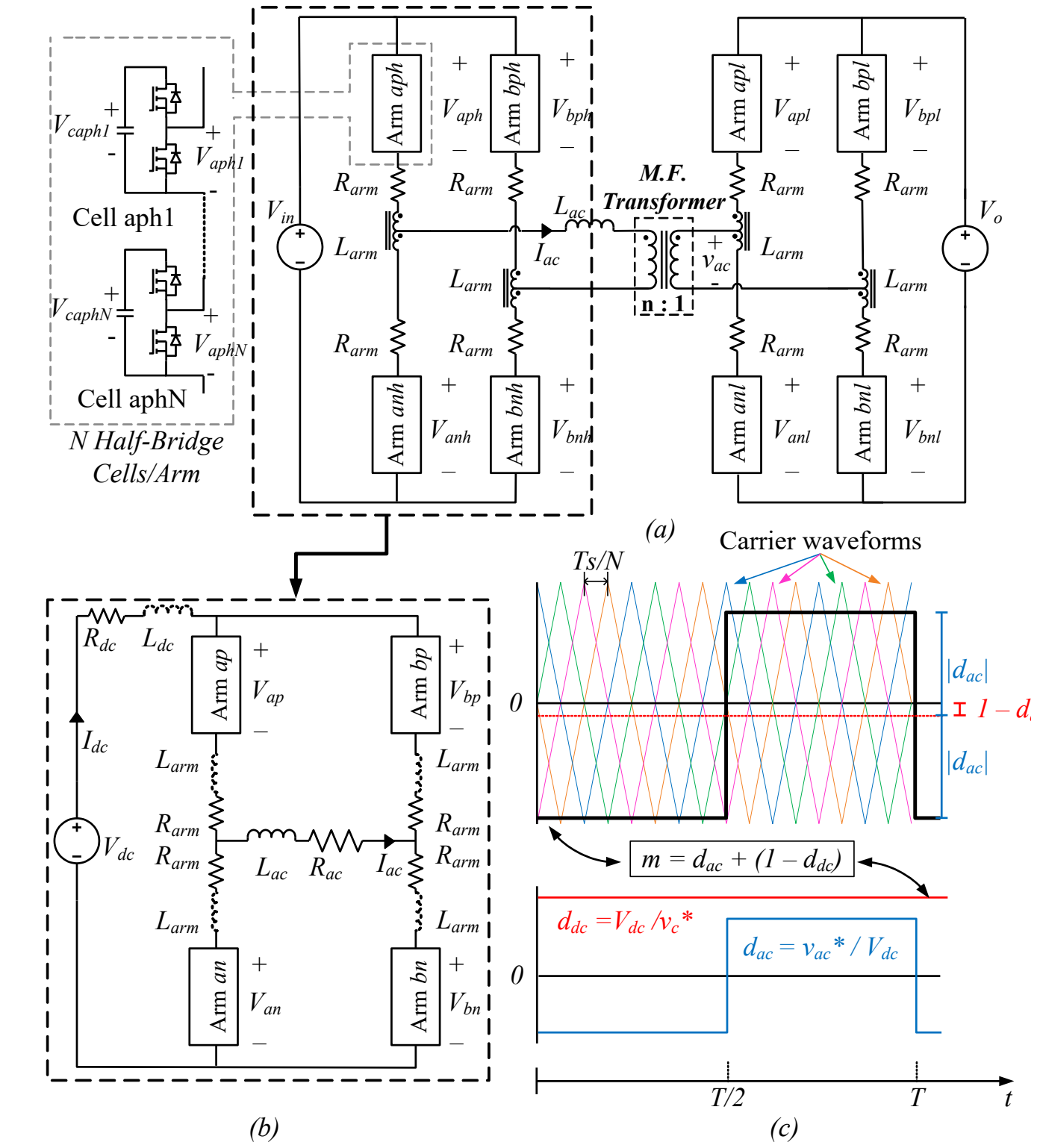

Fig. 1: (a) DC-DC MMC-based SST, (b) one-sided inductor-less simplification with parasitic inductances, and (c) modulation reference and its components.

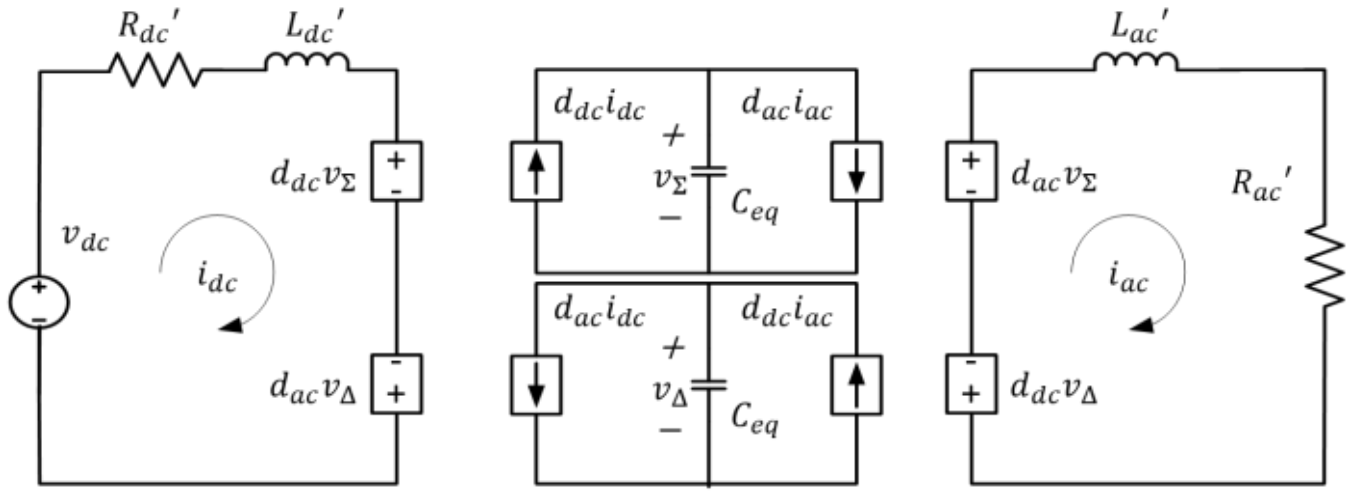

Fig. 2: Equivalent circuit of averaged model from (4).

To derive the averaged model, the state equations (2) are first obtained from the Fig. 1(b) circuit, where $v_{cxy}$ are the arm capacitor voltages and $d_{xy}$ are arm modulation references. A set of definitions are provided in (3). By substituting (3) into (2), the new averaged model is obtained in (4). For readability, the parameters $L_{dc}' = L_{dc} + L_{arm}$, $L_{ac}' = L_{ac} + L_{arm}$, $R_{dc}' = R_{dc} + R_{arm}$, and $R_{ac}' = R_{ac} + R_{arm}$ are defined.

An equivalent circuit based on (4) is shown in Fig. 2, where $C_{eq} = 4C/N$. The variables $v_\Sigma$ and $v_\Delta$ represent the averages of the sum and difference, respectively, of the upper and lower arm capacitor voltages. By defining these variables, it is possible to distinguish the origins of the voltage fluctuations that the equivalent arm capacitances in the MMC experience. The $v_\Sigma$ state variable is associated with the main power transfer path and can be characterized as a DC value superimposed with a 2ω voltage fluctuation. The $v_\Delta$ variable is associated with the MMC circulating energy and experiences a 1ω voltage fluctuation.

To prove the validity of the dynamic model presented in (4), a simulation is performed to compare a circuit-based simulation of the MMC with a simulation of (4) run simultaneously. Fig. 3 is a comparison of the open-loop circuit-based simulation (green, solid) with the open-loop averaged model (blue, dashed). The simulation is performed using an arm inductor-less 8-cell MMC with output frequency $f = 2.5$ *kHz* and carrier frequency $f_{sw} = 17.5$ *kHz*. The currents are normalized to $i_{dc}$* and $|i_{ac}$*$|$, and the voltages are normalized to $v_\Sigma$*. A step change in $|d_{ac}|$ and $d_{dc}$ occurs at $t = 2$ *ms*.

In Fig. 3, the dynamic model simulation and the circuit simulation are in very close agreement both in the steady state and when the step change in control variables is executed. This close agreement is indicative of the validity of the developed dynamic model. Also in Fig. 3, the DC current experiences a very large 2ω ripple. This is caused by the low MMC inductance: significant circulating energy is present because of the lack of inductive damping. This highlights the challenge of

$$\begin{cases} (L_{dc}+L_{arm})\dfrac{d}{dt}\left(\dfrac{i_{ap}+i_{an}+i_{bp}+i_{bn}}{2}\right)=v_{dc}-\dfrac{1}{2}(v_{ap}+v_{an}+v_{bp}+v_{bn})-(R_{dc}+R_{arm})\left(\dfrac{i_{ap}+i_{an}+i_{bp}+i_{bn}}{2}\right) \\ \dfrac{C}{N}\dfrac{d}{dt}\begin{bmatrix} v_{cap}\\ v_{can}\\ v_{cbp}\\ v_{cbn}\end{bmatrix}=\begin{bmatrix} d_{ap} & d_{an} & d_{bp} & d_{bn}\end{bmatrix}\begin{bmatrix} i_{ap}\\ i_{an}\\ i_{bp}\\ i_{bn}\end{bmatrix} \\ (L_{ac}+L_{arm})\dfrac{d}{dt}\left(\dfrac{i_{ap}-i_{an}-i_{bp}+i_{bn}}{2}\right)=-\dfrac{1}{2}(v_{ap}-v_{an}-v_{bp}+v_{bn})-(R_{ac}+R_{arm})\left(\dfrac{i_{ap}-i_{an}-i_{bp}+i_{bn}}{2}\right) \end{cases} \quad (2)$$

$$\begin{cases} d_{ap}=d_{bn}=\dfrac{d_{dc}-d_{ac}}{2} \\ d_{an}=d_{bp}=\dfrac{d_{dc}+d_{ac}}{2}' \end{cases} \begin{cases} i_{ap}=i_{bn}=\dfrac{i_{dc}+i_{ac}}{2} \\ i_{an}=i_{bp}=\dfrac{i_{dc}-i_{ac}}{2}' \end{cases} \begin{cases} v_\Sigma=\dfrac{1}{2}(v_{cap}+v_{can})=\dfrac{1}{2}(v_{cbn}+v_{cbp}) \\ v_\Delta=\dfrac{1}{2}(v_{cap}-v_{can})=\dfrac{1}{2}(v_{cbn}-v_{cbp}) \end{cases} \quad (3)$$

$$\begin{cases} L_{dc}'\dfrac{di_{dc}}{dt}=v_{dc}-d_{dc}v_\Sigma+d_{ac}v_\Delta-R_{dc}'i_{dc} \\ \dfrac{4C}{N}\dfrac{dv_\Sigma}{dt}=d_{dc}i_{dc}-d_{ac}i_{ac} \\ \dfrac{4C}{N}\dfrac{dv_\Delta}{dt}=d_{dc}i_{ac}-d_{ac}i_{dc} \\ L_{ac}'\dfrac{di_{ac}}{dt}=-d_{dc}v_\Delta+d_{ac}v_\Sigma-R_{ac}'i_{ac} \end{cases} \quad (4)$$

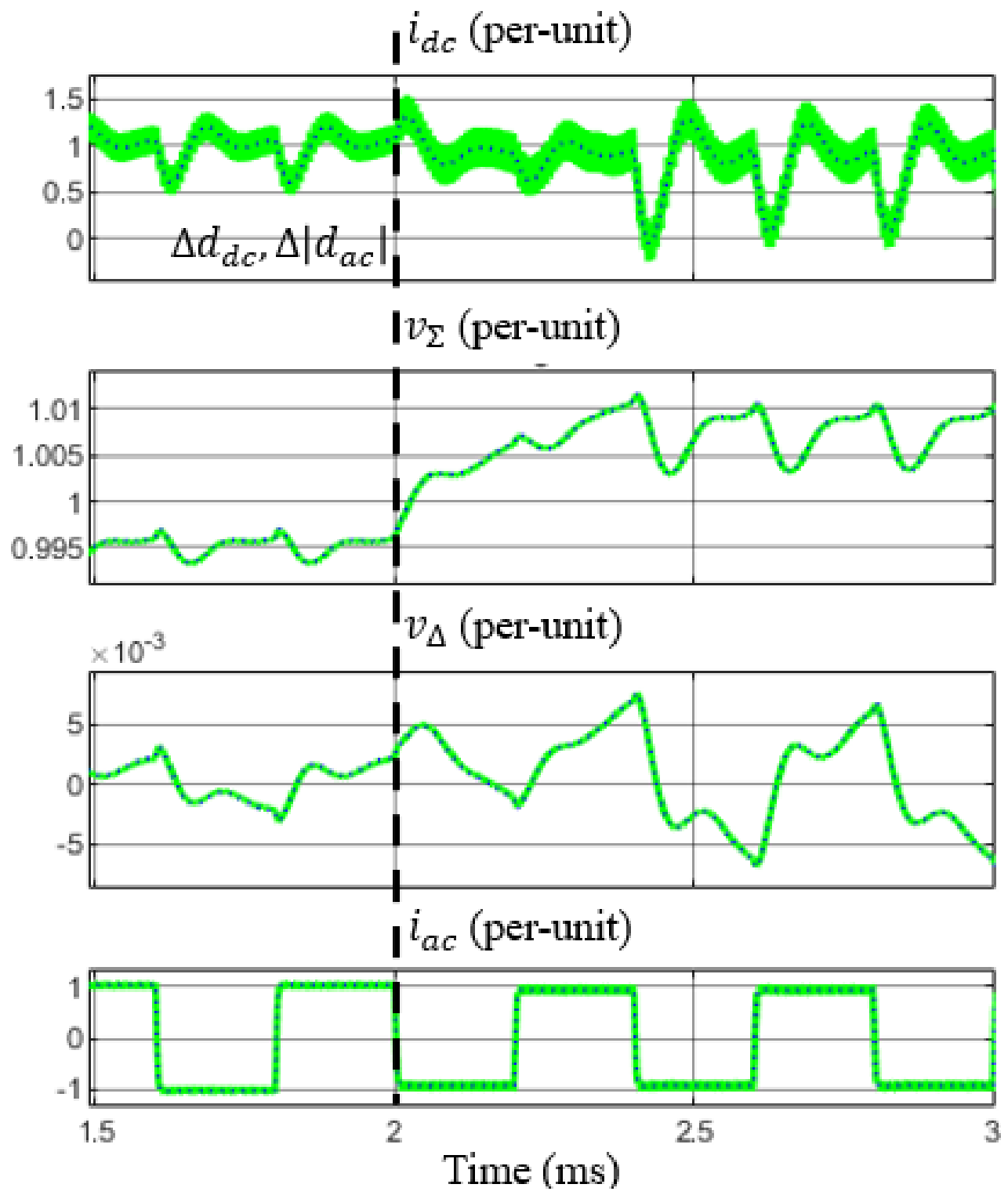

Fig. 3: Open-loop comparison of a circuit-based simulation (green, solid) and the averaged model (blue, dashed). A step change in $|d_{ac}|$ and $d_{dc}$ occurs at *t = 2 ms.*

designing a controller for a low-inertia MMC, as the designed controller must damp this energy while performing the other control functions.

## III. Modulated Model Predictive Controller

In this section, a deadbeat-type MPC with low computational burden is designed based on the developed averaged model. The controller is shown in Fig. 4, where the power flow control is implemented in the predictive controller, and the circulating energy control is implemented both in the predictive controller and in the steady state compensator. The controller samples the DC voltage, capacitor voltages, DC current, and AC current from the converter, and receives references for the AC current and $v_\Sigma$. After determining the DC current reference through the controller $G_c^{vi}$ and predicting $v_\Delta$, the controller calculates $d_{ac}(k+1)$, $d_{dc}(k+1)$, and the compensator output, and then generates the modulation waveform. A cell-swapping algorithm is implemented to balance the cell capacitors within each arm, which is described in the appendix. The AC current reference is determined from the power demand, the capacitor voltage reference is based on the DC input voltage, and the DC current reference is generated by a model-based PI controller, $G_c^{vi}$, whose design is explained in section B and is based on the method from [26]. This simple implementation is far less complicated than typical MPCs and is enabled in part by the hybrid nature of the controller.

The rest of this section is organized as follows: the predictive controller is described in section A; the model-based PI controller for generating the DC current reference is described in section B; and the model-based feed-forward compensator is described in section C.

### A. Predictive Controller

The predictive controller performs the power flow control and part of the circulating energy control. To design the controller, the DC and AC current state equations from (4) are discretized using trapezoidal integration as shown in (5)-(6), where $T_s$ is the sampling period. This controller will minimize the error between the sensed values and references of $i_{dc}$ and $i_{ac}$, where the $i_{ac}$ reference is generated based on the power demand

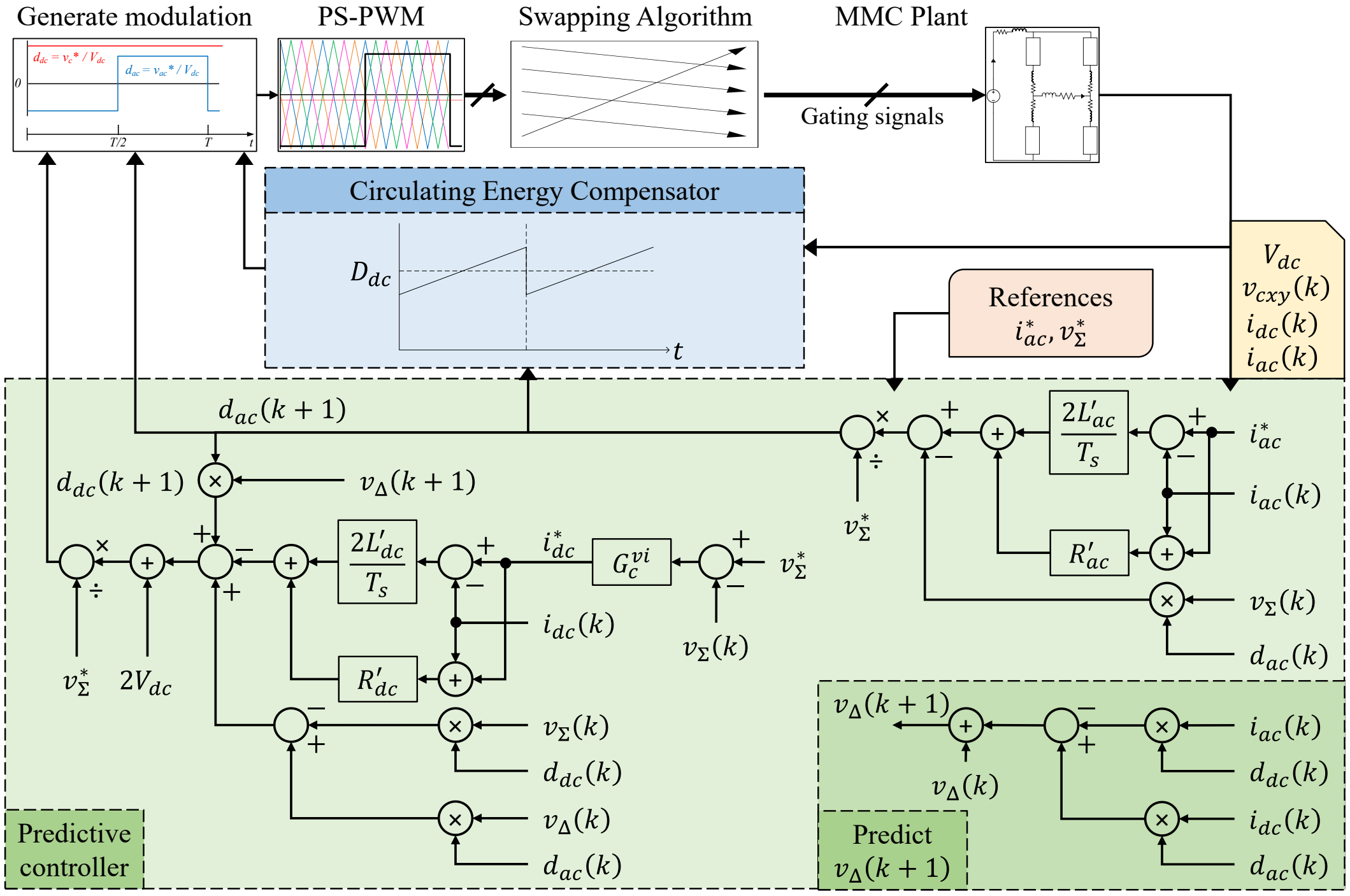

Fig. 4: Control block diagram of proposed modulated model predictive controller with circulating energy controller.

$$\begin{aligned}\left(i_{ac}(k+1)-i_{ac}(k)\right) \\ &= \frac{T_s}{2L'_{ac}}\left[-\left(d_{dc}(k+1)v_\Delta(k+1)+d_{dc}(k)v_\Delta(k)\right)+\left(d_{ac}(k+1)v_\Sigma(k+1)+d_{ac}(k)v_\Sigma(k)\right)\right. \\ &\left.-R'_{ac}\left(i_{ac}(k+1)+i_{ac}(k)\right)\right]\end{aligned} \tag{5}$$

$$\begin{aligned}\left(i_{dc}(k+1)-i_{dc}(k)\right) \\ &= \frac{T_s}{2L'_{dc}}\left[V_{dc}-\left(d_{dc}(k+1)v_\Sigma(k+1)+d_{dc}(k)v_\Sigma(k)\right)\right. \\ &\left.+\left(d_{ac}(k+1)v_\Delta(k+1)+d_{ac}(k)v_\Delta(k)\right)-R'_{dc}\left(i_{dc}(k+1)+i_{dc}(k)\right)\right]\end{aligned} \tag{6}$$

$$d_{ac}(k+1)=\left(2L'_{ac}\left(i^*_{ac}-i_{ac}(k)\right)/T_s+R'_{ac}\left(i^*_{ac}+i_{ac}(k)\right)-d_{ac}(k)v_\Sigma(k)\right)/v^*_\Sigma \tag{7}$$

$$\begin{aligned}d_{dc}(k+1)=&\left(2V_{dc}+d_{ac}(k+1)v_\Delta(k+1)-2L'_{dc}\left(i^*_{dc}-i_{dc}(k)\right)/T_s-R'_{dc}\left(i^*_{dc}+i_{dc}(k)\right)-d_{dc}(k)v_\Sigma(k)\right. \\ &\left.+d_{ac}(k)v_\Delta(k)\right)/v^*_\Sigma\end{aligned} \tag{8}$$

$$v_\Delta(k+1)=\frac{NT_s}{4C}\left(d_{dc}(k)i_{ac}(k)-d_{ac}(k)i_{dc}(k)\right)+v_\Delta(k) \tag{9}$$

$$G_{vi}=\frac{\hat{v}_\Sigma}{\hat{\imath}_{dc}}=\frac{V_\Sigma D_{dc}R_{ac}-I_{dc}R_{dc}R_{ac}}{I_{dc}D_{dc}R_{ac}+V_\Sigma|D_{ac}|^2}\times\frac{s^2\frac{I_{dc}L_{dc}L_{ac}}{I_{dc}R_{dc}R_{ac}-V_\Sigma D_{dc}R_{ac}}+s\left(\frac{I_{dc}(R_{ac}L_{dc}+R_{dc}L_{ac})-V_\Sigma D_{dc}L_{ac}}{I_{dc}R_{dc}R_{ac}-V_\Sigma D_{dc}R_{ac}}\right)+1}{s^2\frac{V_\Sigma L_{ac}C}{I_{dc}D_{dc}R_{ac}+V_\Sigma|D_{ac}|^2}+s\left(\frac{V_\Sigma R_{ac}C+I_{dc}D_{dc}L_{ac}}{I_{dc}D_{dc}R_{ac}+V_\Sigma|D_{ac}|^2}\right)+1} \tag{10}$$

and the $i_{dc}$ reference is generated based on the capacitor voltage control described in section B. To this end, $i_{dc}(k+1)$, $i_{ac}(k+1)$, and $v_\Sigma(k+1)$ are replaced with their references: $i_{dc}(k+1) = i_{dc}^*$, $i_{ac}(k+1) = i_{ac}^*$, and $v_\Sigma(k+1) = v_\Sigma^*$. This assumption states that the controller is fast and accurate enough to reach the desired value within one time step, which permits a closed-form solution.

Equation (5) is rearranged to solve for $d_{ac}(k+1)$ in (7) to perform the power flow control, while (6) is rearranged to solve for $d_{dc}(k+1)$ in (8) to perform the circulating energy control. In obtaining (7) from (5), the $v_\Delta$ terms are discarded. These terms represent the converter's circulating energy and are negligible compared to other terms in (7). Furthermore, because (7) corresponds to the power flow control, the circulating energy information can be neglected with the assumption that it will be well-controlled by the circulating energy controller in (8), similar to the *dq*-based approach in traditional MMC power flow control. The predicted value of $v_\Delta$, $v_\Delta(k+1)$, is needed in (8) and is obtained from a Euler discretization, shown in (9).

The predictive controller is implemented sequentially: $d_{ac}(k+1)$ is first calculated to achieve the power flow control, and then $d_{dc}(k+1)$ is calculated as part of the circulating energy control. The calculation of $d_{dc}(k+1)$ is analogous to the inner current loop of a traditional MMC circulating energy controller such as in [11] and is related to the arm voltage control. The method for obtaining the circulating current reference ($i_{dc}^*$ in this paper) is described in section B. An additional controller component is introduced in section C to reduce the significant steady-state circulating energy in this low-inertia converter.

## B. Circulating Current Reference Generation

MMC controllers in both traditional and MPC-based approaches typically use a cascaded structure for the circulating energy control, with an outer voltage loop and an inner current loop. A similar strategy is adopted in this paper, but unlike other approaches, the outer loop controller in this paper is model-based and does not require parametric tuning. First, the circulating energy controller is shown in Fig. 5(a), where the predictive inner loop controller is nonlinear. To design the outer loop controller, the approximation $i_{dc}^* \approx i_{dc}$ is made to simplify the structure to Fig. 5(b), which is linear time-invariant, and where $G_{vi}$ is the transfer function from $i_{dc}$ to $v_\Sigma$. $G_{vi}$ is described in (10) and derived according to [26]. A controller $G_c^{vi}(s) = K_p\left(1+\frac{K_i}{s}\right)$ to generate $i_{dc}^*$ can be designed by examining the frequency response of this transfer function without requiring trial-and-error tuning.

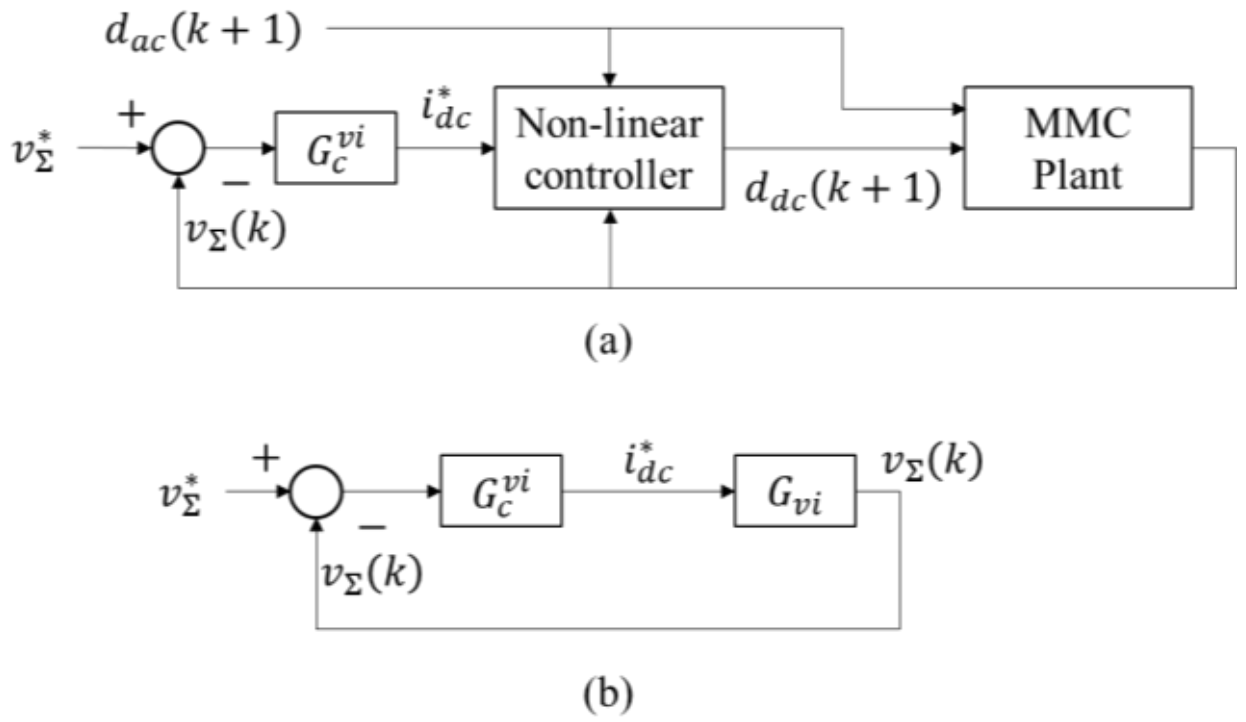

Fig. 5: (a) Controller structure showing the outer loop and nonlinear inner loop; (b) simplification of the circulating energy controller to a linear time-invariant system for the design of the outer loop controller.

## C. Circulating Energy Compensator

In the low-inertia MMC, control of the steady-state circulating energy is critical, especially because the square-wave modulation will cause significant resonance every half-period. This problem could be solved by increasing the sampling speed of the predictive controller in (8) and carefully defining a differential voltage reference $v_\Delta^*$ as part of the circulating current reference. In this section, an alternative approach is examined for reducing the steady-state circulating energy based on a fast model-based compensator.

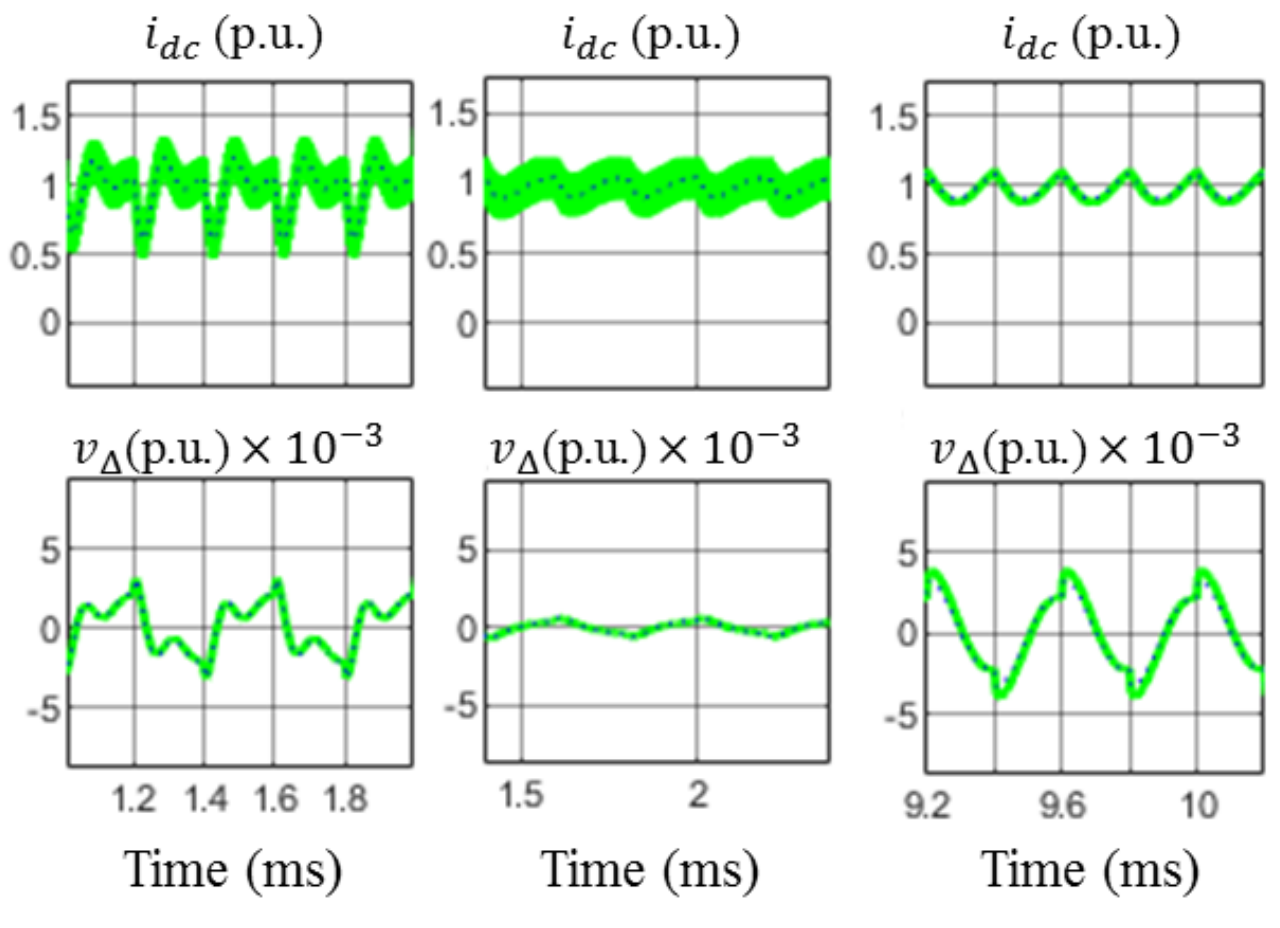

Fig. 6: Open-loop simulation of circuit-based model (green, solid) and dynamic model (blue, dashed) showing $i_{dc}$ and $v_\Delta$ using differently sized passive components: (a) small inductors and small capacitors; (b) small inductors and large capacitors; and (c) large inductors and small capacitors.

To better understand how circulating energy manifests in this converter, it is useful to observe the $i_{dc}$ and $v_\Delta$ state variables. Three open-loop simulation cases are shown in Fig. 6, differing only in passive component sizing: (a) represents the inductor-less case with small cell capacitors; (b) represents the inductor-less case with large cell capacitors; and (c) uses large inductors and small cell capacitors. The $i_{dc}$ waveform is normalized to $i_{dc}$*, and the $v_\Delta$ waveform is normalized to $v_\Sigma$*. The effect of the circulating energy on $i_{dc}$ is obvious: with small passive components, the open-loop $i_{dc}$ experiences a large double-frequency ripple due to the circulating energy, but not in the other two cases.

We notice that in Fig. 6(c), the case with large arm inductors, the arm differential voltage $v_\Delta$ has the same magnitude as in the inductor-less case, Fig. 6(a). Unlike the inductor-less case, however, the arm differential voltage appears as a roughly triangular waveform. From this heuristic observation, we propose the following: if the arm differential voltage $v_\Delta$ can be controlled to assume a triangular waveform, then the MMC will behave as if it were subject to a large arm inductance – i.e. a virtual inductance is designed – and the circulating energy will be correspondingly reduced.

For clarity, the controller proposed in this section will henceforth be referred to as the "steady state compensator." The derivation of this proposed compensator begins with the dynamic model. First, the $i_{dc}$ state equation is reproduced in (11) with the resistances neglected. It is equated to zero, asserting that $i_{dc}$ is constant and thus contains no distortion caused by circulating energy. This is the desired result of the controller, and the controller is derived backwards to meet this condition. Because the proposed compensator aims to reduce the circulating energy in the MMC, it will use the $d_{dc}$ control variable. A graphical representation of how this compensator's output interacts with the output of the predictive controller from section A is shown in Fig. 7.

Equation (11) is solved for $d_{dc}$ in (12), which separates $d_{dc}$ into the component generated by the predictive controller and the component generated by the steady-state compensator according to Fig. 7.

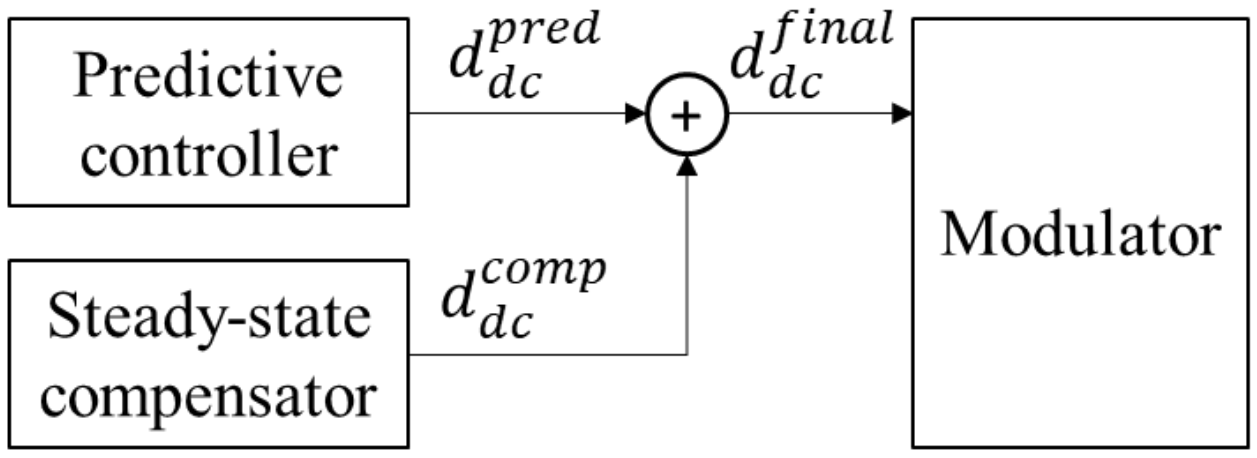

Fig. 7: Interaction of the predictive controller (section A) and the steady-state compensator (section C) $d_{dc}$ outputs. The compensator is simple and can operate much faster than the predictive controller.

$$\begin{cases} L'_{dc}\dfrac{di_{dc}}{dt} = v_{dc} - d_{dc}v_\Sigma + d_{ac}v_\Delta - R'_{dc}i_{dc} \\ 0 = v_{dc} - d_{dc}v_\Sigma + d_{ac}v_\Delta \end{cases} \tag{11}$$

$$\begin{cases} d_{dc}^{final}(t) = \dfrac{v_{dc}}{v_\Sigma(t)} + \dfrac{d_{ac}(t)v_\Delta(t)}{v_\Sigma(t)} \\ d_{dc}^{final}(t) = d_{dc}^{pred} + d_{dc}^{comp} \end{cases} \tag{12}$$

From (12), the expression for the compensator output requires knowledge of $v_\Sigma(t)$ and $v_\Delta(t)$, while $d_{ac}(t)$ is already known from the output of the predictive controller from section A. If the system is well-controlled, $v_\Sigma(t)$ can be approximated as $v_\Sigma^*$, leaving only $v_\Delta(t)$. As noted earlier in this section, the basis of this compensator lies in shaping $v_\Delta$ into a triangular waveform, so $v_\Delta(t)$ can be expressed as the desired triangular function. To obtain an expression for $v_\Delta$, its amplitude must be known. To do so, the $v_\Delta$ state equation from (4) is rewritten as a difference expression in (13), where $\Delta T$ is one half of the AC period, $I_{dc}$ is the average steady-state DC current, and $I_{ac}$ is the RMS AC current. Note that a triangular waveform has a constant slope in each half of its period, which is why (13) is a valid expression. Then, using simple algebra for describing a triangle wave and simplifying the resultant expression, the compensator's output may be written as in (14), where $t$ is periodic in $T$/2 and increases linearly. The controller output is demonstrated in Fig. 8. Note that because only the average or RMS values of $I_{dc}$ and $I_{ac}$ are needed for this compensator, it can operate much faster than the predictive controller.

$$\frac{4C}{N}\frac{\Delta V_\Delta}{\Delta T} = d_{dc}^{pred}I_{ac} - \left|d_{ac}^{pred}\right|I_{dc} \tag{13}$$

$$d_{dc}^{comp} = \frac{\left|d_{ac}^{pred}\right|}{v_\Sigma^*}\left(\frac{2}{T}t - \frac{1}{2}\right)\frac{T}{2}\left(\frac{d_{dc}^{pred}I_{ac} - \left|d_{ac}^{pred}\right|I_{dc}}{4C/N}\right) \tag{14}$$

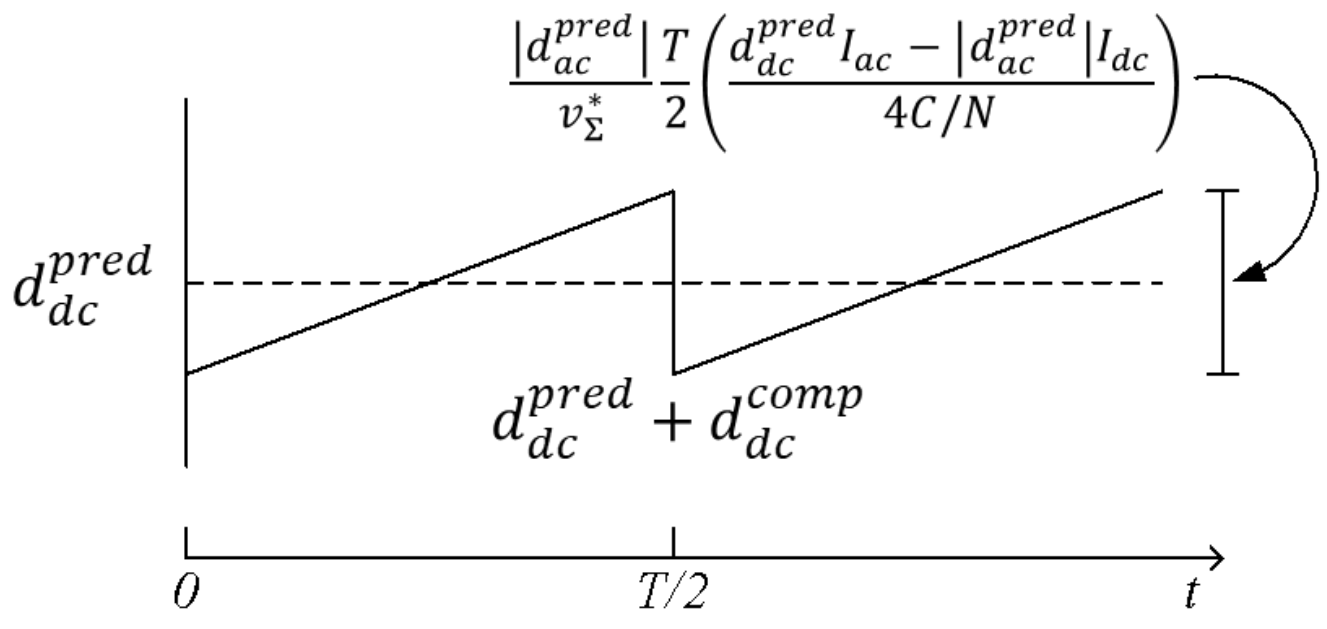

Fig. 8: Output of the steady state circulating energy controller.

This feed-forward controller is comparable to an indirect modulation control, or a "virtual inductance" control. In indirect modulation, the MMC reference voltages are normalized using the sensed or predicted values of the capacitor voltages, thus incorporating the circulating energy information in the modulation and thereby reducing it [5]. The novelty of the proposed circulating energy controller, besides its application to a square-wave modulation, is its simplicity, its speed, its model-based nature, and its equivalent description as a virtual inductance.

## IV. Simulation Case Study

To validate the proposed control system, an MVDC MMC case study is performed in a Simulink environment. The simulation of one MMC is presented, representing one side of the full DC-DC SST. The simulation parameters are listed in Table 1, corresponding to the definitions shown in Fig. 1(b). The outer loop controller achieves a 3 kHz bandwidth, while the predictive inner loop bandwidth is 50 kHz.

TABLE I
MMC Simulation Parameters

| Variable | Per Unit Values | Actual value |
|---|---|---|
| $V_{dc}$ | - | 6 kV |
| $P_{rated}$ | - | 500 kW |
| Cells per arm | - | 8 |
| $f_{ac}$ | - | 2.5 kHz |
| $f_{sw}$ | - | 40 kHz |
| C | 1.13 (arm) | 100 μF (cell) |
| $R_{ac}$ | - | 65 Ω |
| $R_{dc}$ | 0.5 % | 0.325 Ω |
| $R_{arm}$ | 0.1 % | 0.065 Ω |
| $L_{ac}$ | 3.0 % | 124 μH |
| $L_{dc}$ | 0.2 % | 8.28 μH |
| $L_{arm}$ | 0.1 % | 4.14 μH |
| $f_{sample}$ | - | 50 kHz |
| $K_p$ | - | 1 |
| $K_i$ | - | 2π×1000 Hz |

Simulation results are presented in Fig. 9. Fig. 9(a)-(c) each consist of five windows. The topmost window shows the AC current (blue) compared to its reference (black, bold). The second window shows $v_\Sigma$ (blue) compared to its reference (black, bold). The third window shows the DC current, the fourth shows the generated $d_{dc}$ control variable, and the fifth shows the generated $d_{ac}$ control variable. The currents and voltages in Fig. 9 are all from the circuit simulation, not from the dynamic model. Fig. 9(a)-(c) show a time span of 3 ms.

Fig. 9(a) shows the converter initially in open-loop operation. Prior to the controller activation, the DC current suffers from a significant double-frequency ripple component, causing it to fluctuate from 45 A to 105 A. At *t = 4 ms*, the control system is activated, and the DC current becomes well-regulated as the circulating energy is controlled. Though a ripple component is still present, it is significantly reduced by almost 60%. Additionally, the small error in $v_\Sigma$ is corrected.

Fig. 9(b) shows the converter initially in the steady state at rated power with the control system active. At *t = 10 ms*, there is a step change in the power reference. The power flow control responds accordingly by reducing the magnitude of $d_{ac}$, and the circulating-energy control slightly increases the magnitude of $d_{dc}$. There is a very small fluctuation in $v_\Sigma$ at the transition, which is quickly corrected, and within 3 AC periods, the DC current becomes tightly regulated once again.

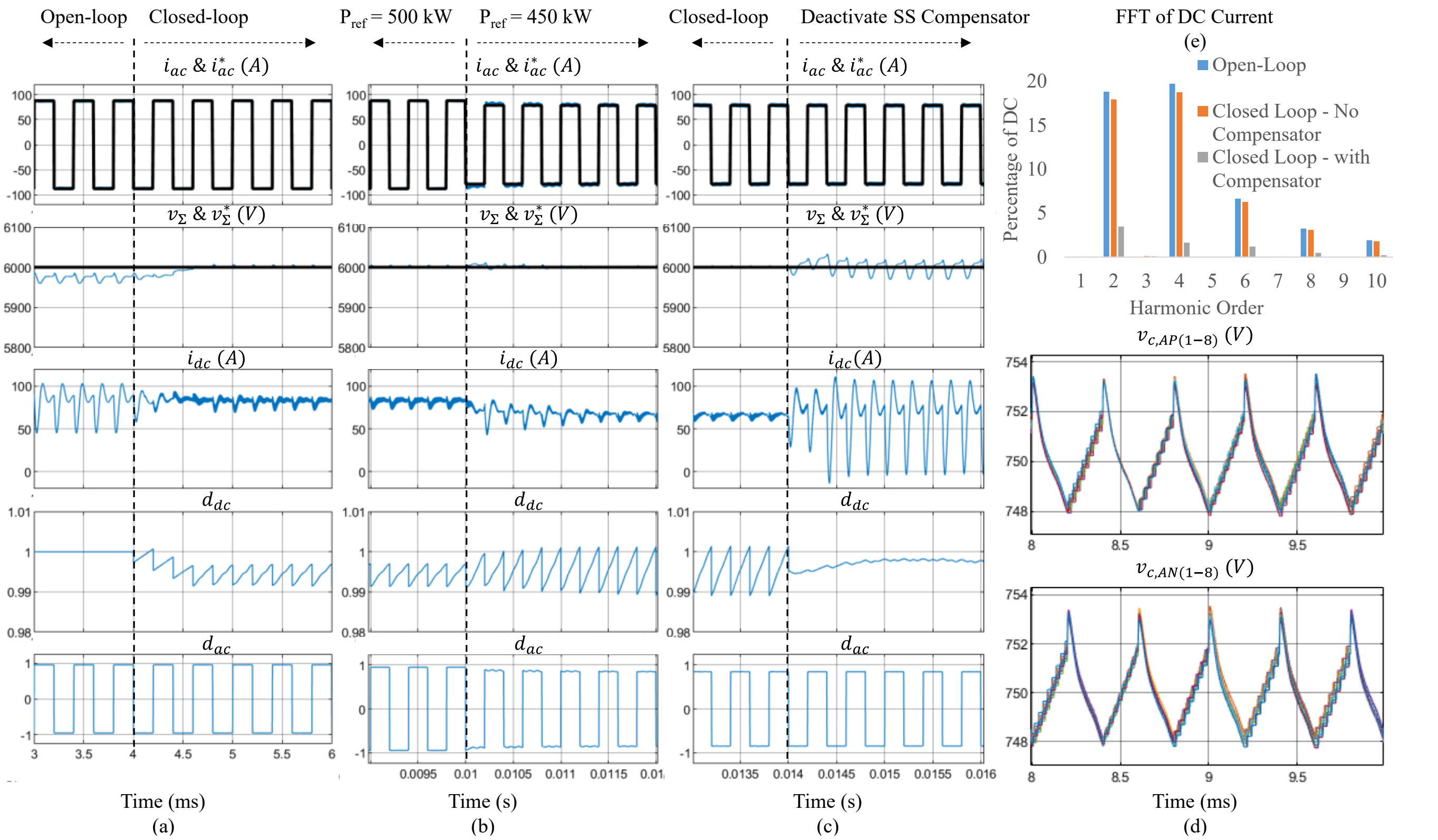

Fig. 9: Performance of proposed deadbeat-style modulated MPC: (a) activation from open loop; (b) transient response to power reference step change; (c) deactivation of the circulating energy compensator; (d) closed-loop capacitor voltages in phase A; and (e) FFT of the DC current under different conditions: open loop (blue, left), closed loop w/o compensator (orange, middle), and closed loop w/ compensator (gray, right).

Fig. 9(c) shows the converter initially in the steady state with the control system active. At *t = 14 ms*, the steady state circulating energy compensator is deactivated, leaving just the predictive controller active. Immediately, the DC current waveform degrades significantly and experiences a very large double-frequency ripple, indicating that the compensator was having a significant impact on the circulating energy suppression. With the increase in the DC current ripple, $v_\Sigma$ also begins to experience a small double-frequency ripple.

Fig. 9(d) shows the capacitor voltages for the MMC arms in one phase leg and demonstrates that the open-loop cell-swapping algorithm is performing the individual cell balancing task well. Fig. 9(e) further highlights the performance of the controller through an FFT analysis of the DC current under three different conditions: open loop (blue, left), closed loop without the compensator (orange, middle), and closed loop with the compensator (gray, right). Without control, the arm inductor-less square-wave MMC experiences large $2^{nd}$-order and $4^{th}$-order harmonic components of almost 20% the DC value. The predictive control is not designed to control the steady state circulating energy, so the harmonics are relatively unchanged when it is active. On the other hand, the low frequency harmonics of the DC current are significantly mitigated from almost 20% to less than 4% by the steady state compensator, validating its usefulness and its powerful capabilities. The simulation results have demonstrated that the proposed MPC is a viable strategy for controlling the arm inductor-less MMC, and that the novel steady-state compensator can significantly mitigate the circulating energy in the inductor-less MMC.

For comparison, the simulation is repeated with the controller proposed in [17], with comparison results shown in Fig. 10(a)-(b). In both Fig. 10(a) and (b), the power reference is stepped from 500 kW to 450 kW at $t$ = 14 ms. Fig. 10(a) is a simulation of the controller proposed in this paper, while Fig. 10(b) shows a simulation using the MPC in [17] under the same conditions. For Fig. 10(b), the circulating current reference is generated similarly to [17], although the PI parameters are not an exact match due to the different AC frequency and modulation.

From Fig. 10(a)-(b), when the same sampling speed is used for both controllers, the controller from [17] fares significantly worse than the proposed MPC. Fig. 10(b) suffers from very large low-frequency DC current harmonics, especially after the power reference step change, and suffers from a larger voltage overshoot during the transient. The controller from [17] uses indirect modulation to reduce the MMC circulating energy, but it requires a much higher sampling speed than the method proposed in this paper. To further highlight this, another simulation was performed with the MPC from [17], this time with a much higher sampling frequency of 300 kHz, six times greater than the sampling frequency used in Fig. 9 and Fig. 10. The results are presented in Table II. With this high sampling speed, the performance of the controller from [17] is comparable to the proposed controller. However, it still suffers a greater voltage overshoot, and the current control is also inferior. Since the controller proposed in this paper achieves a better result than the controller from [17] despite operating with a 6 times smaller sampling frequency, the superiority of the proposed controller is established.

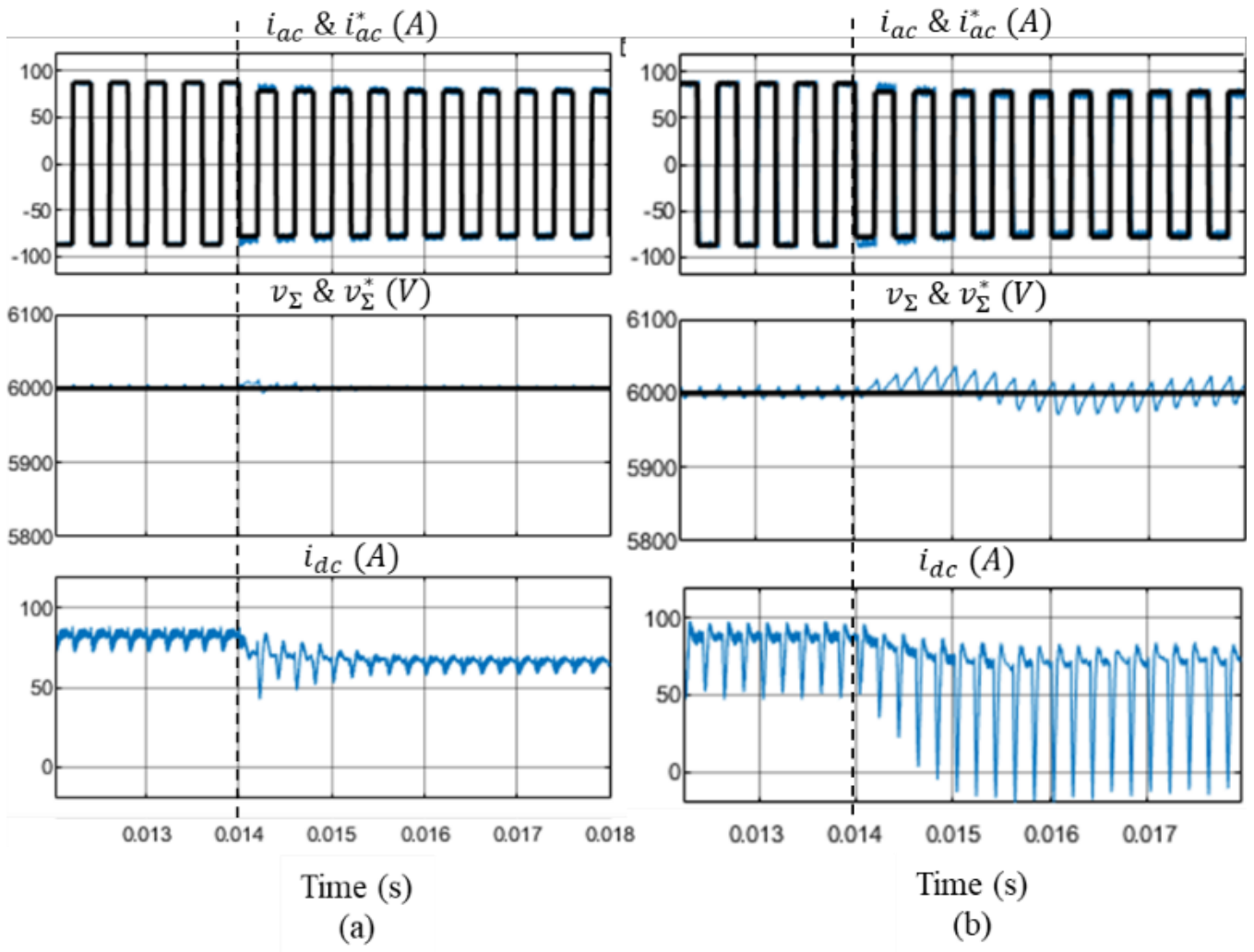

Fig. 10: Simulation comparison of (a) the proposed controller and (b) the controller from [17], both with a sampling frequency of 50 kHz. At t = 14 ms, a power reference step change occurs.

TABLE II
Controller Comparison

| Value | Proposed Controller | Controller in [17] | Controller in [17] – fast |
|---|---|---|---|
| fsample | 50 kHz | 50 kHz | 300 kHz |
| $v_\Sigma$ overshoot | ~0.2% | ~0.7% | ~0.7% |
| $i_{dc}$ 2ω FFT | 3.5 % of DC | 11.9 % of DC | 4.5 % of DC |
| $i_{dc}$ 4ω FFT | 1.6 % of DC | 12.3 % of DC | 2.3 % of DC |

## V. Experimental Verification

A downscaled MMC hardware prototype is built in the laboratory and shown in Fig. 11, with its hardware parameters listed in Table III. The single-phase MMC consists of 16 cells arranged into 4 arms. No arm inductors are used in the MMC prototype; the DC inductance listed in Table III is an estimate of the parasitic inductance in the circuit.

The control system was implemented in a MyWay PE-Expert 4, which features a TI C6657 DSP and a Xilinx XC7K70T FPGA, the latter of which is operated on a 100 MHz clock. The DSP is operated with a 100 kHz sampling speed. The predictive controller is implemented in C code on the DSP, while the carrier generation, PWM, and cell-swapping algorithm are implemented in VHDL on the FPGA.

Experimental waveforms are presented in Fig. 12(a)-(c). From the top, Fig. 12(a) and (c) show four capacitor voltages from one arm, then the AC voltage (yellow), the AC current (green), and the DC current (magenta). Fig. 12(b) shows eight capacitor voltages from one phase leg. The voltage and current divisions are listed in the figures.

Fig. 12(a) begins with the system in the open loop, with the DC current significantly distorted by a double-frequency ripple. At the time indicated by the dashed line, the system enters the closed loop at full power, with the DC current immediately becoming well-regulated as the controller reduces the circulating energy in the inductor-less converter.

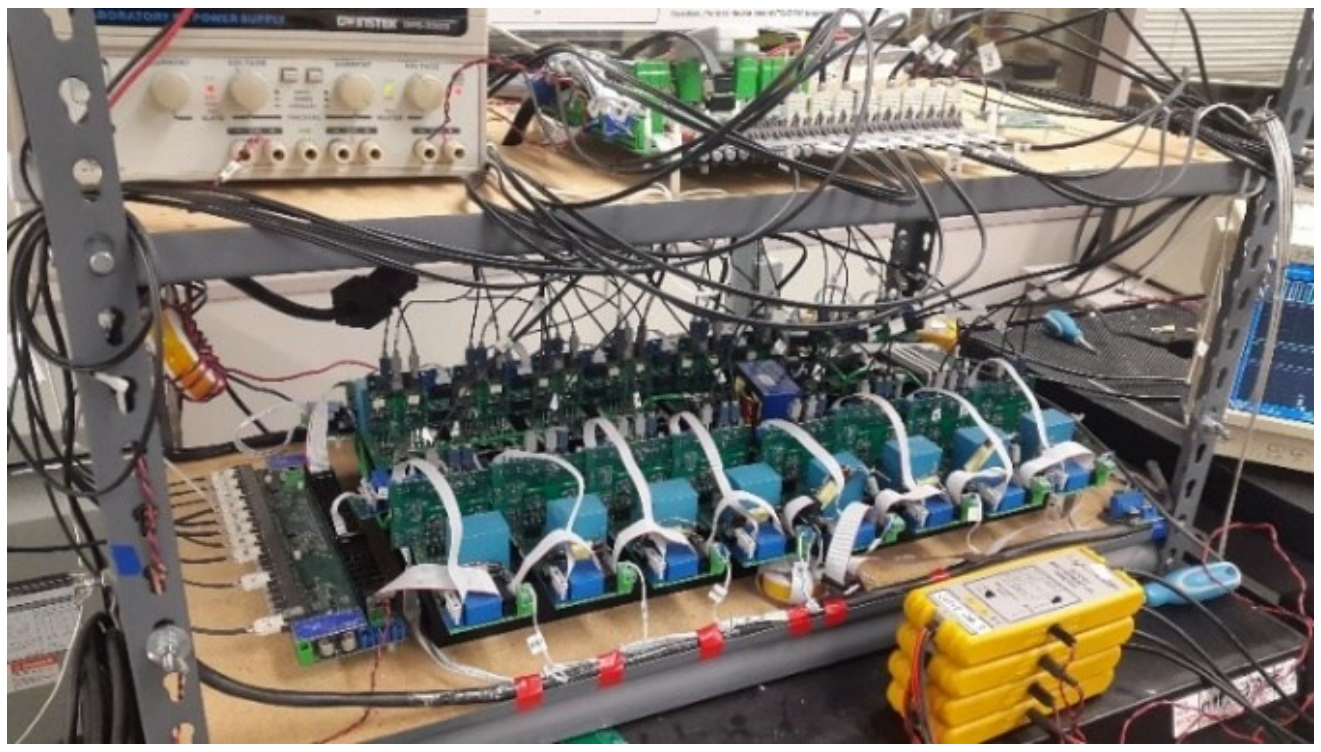

Fig. 11: Laboratory MMC hardware.

TABLE III
MMC Hardware Parameters

| Variable | Per Unit Values | Actual value |
|---|---|---|
| Vdc | - | 250 V |
| P | - | 2.5 kW |
| Cells per arm | - | 4 |
| fac | - | 5 kHz |
| fsw | - | 40 kHz |
| C | 0.31 (arm) | 20 µF (cell) |
| Rac | - | 16 Ω |
| Rdc | ~1.25 % | ~0.2 Ω |
| Lac | 8.7 % | 29 µH |
| Ldc | ~0.6 % | ~2 µH |
| fsample | - | 100 kHz |

Fig. 12(c) shows the system's transient response. The converter initially operates at 1 kW. At the time indicated by the dashed line, the power reference is suddenly changed to 2.5 kW, a 150% increase. The controller quickly reacts to increase both the AC current and DC current to meet the power demand and maintain the desired capacitor voltage. It is important to note that the large ripple of $i_{dc}$ in Fig. 12(c) is switching current ripple due to the low power rating and small cell number of the downscaled testbed. For a real high power MMC with a larger number of cells, this large switching ripple can be reduced significantly, as shown in the simulation results of Fig. 9. This large switching frequency ripple cannot be suppressed by control.

Fig. 12(b) offers a closer view of the capacitor voltages during this power reference step change. At the time that the step change in the reference occurs, the capacitor voltages do not experience any appreciable voltage overshoot transient. Fig. 12(d) shows a matching simulation of a recreation of the experimental setup and controller subject to the same conditions. As seen by comparing Fig. 12(c) and Fig. 12(d), the expected behavior of the controller matches the experimental waveforms.

## VI. CONCLUSION

In this paper, a modulated MPC with low computational burden was designed based on a derived comprehensive dynamic model of the MMC for DC-SST application. The MPC performs power flow control and circulating energy control for a low-inertia MMC, whose control is challenging because of its arm inductor-less nature. Using the derived model and its explicit definition of the circulating energy, a fast and highly effective model-based steady-state circulating energy compensator was derived, which operates as a virtual inductance to mitigate low order harmonics of the DC current caused by the increased circulating energy in the arm inductor-less MMC. The proposed controller was verified in a simulation case study and a downscaled MMC hardware prototype, and was also compared to a similar published controller. Both simulation results and experimental results have demonstrated that the proposed control system is a viable, competitive solution to control a low-inertia MMC, and that it can achieve its control objectives with a relatively small sampling speed.

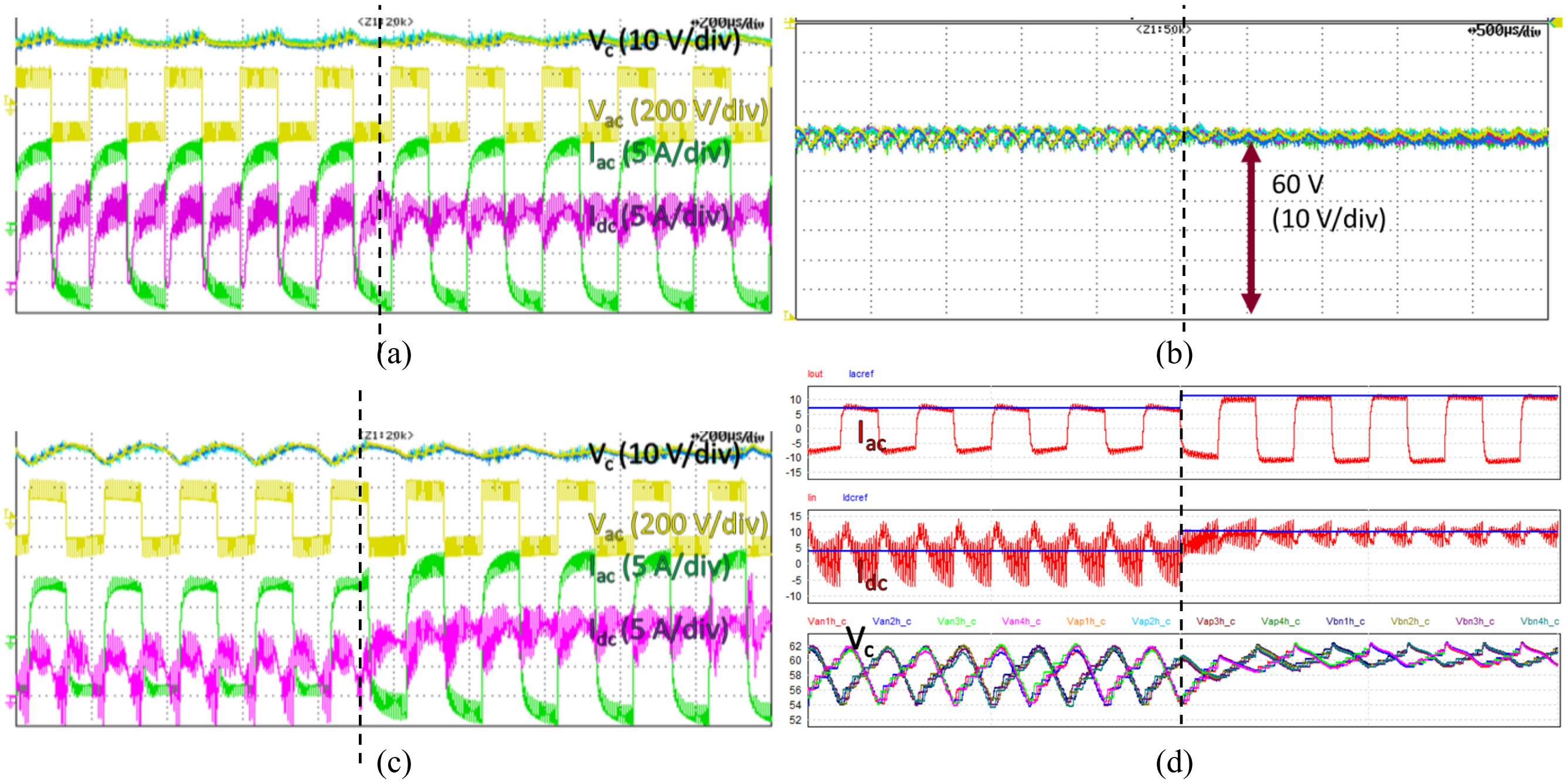

Fig. 12: Experimental waveforms of a downscaled MMC: (a) arm capacitor voltages, AC voltage (yellow), AC current (green), and DC current (magenta) during transition from open loop to closed loop; (b) capacitor voltages in one phase leg during 150% power reference step change; (c) arm capacitor voltages, AC voltage (yellow), AC current (green), and DC current (magenta) during 150% power reference step change; and (d) matching simulation of AC current and reference, DC current and reference, and leg capacitor voltages during 150% power reference step change.

## VII. Appendix

Table A-1 provides a per-unit comparison of the sizes of passive components in selected PI-based and MPC-based MMCs reported so far. The per-unit values were derived from experimental setups where available, and from simulation parameters where not. The last column shows the controller sampling frequencies. Note that the sampling frequency of the proposed MPC is high because its AC frequency is 5 kHz, whereas most other cases operate at 50 Hz or 60 Hz.

Table A-1
MMC Passive Component Size Comparison

| Ref. | Controller | Arm Cap. (p.u.) | Arm Ind. (%) | Sampling Frequency (Hz) |
|---|---|---|---|---|
| [11] | PI | 1.12 | 3.49 | 16 000 |
| [12] | PI | 6.75 | 9.92 | 12 000 |
| [13] | PI | 0.56 | 10.32 | 4 000 |
| [14] | PI | 0.17 | 7.99 | UA |
| [19] | MPC OSS | 0.72 | 1.80 | 1 000 |
| [22] | MPC OSS | 0.83 | 19.63 | 10 000 |
| [24] | MPC OSS | 11.90 | 5.24 | 10 000 |
| [20] | MPC OVL | 1.92 | 4.67 | 4 150 |
| [21] | MPC OVL | 0.86 | 5.20 | 1 500 |
| [23] | MPC OVL | 3.44 | 11.52 | 10 000 |
| [25] | MPC OVL | 0.75 | 7.35 | 20 000 |
| [17] | Deadbeat | 0.64 | 22.2 | 6 000 |
| [18] | Deadbeat | 0.05 | 108.8 | 6 650 |
| Proposed | | 0.31 | ~0.6 | 100 000 |

The cell-swapping algorithm used to balance the capacitor voltages is described in Table A-2 for the experimental setup. The algorithm works in an open-loop manner and can be described as a state machine, where the number of states is equal to twice the number of cells. The entries correspond to MMC cells, where the first letter denotes whether they are in leg A or B, the second letter whether they are in the positive or negative arm, and the number indicates the cell number within that arm. For readability, only the cells in leg A are shown; for leg B, the cells correspond to their opposite arm, i.e. cell BN1 is triggered in the same way as cell AP1, etc. The 'order' column lists the order in which the cells are triggered. For example, in state 0 (S0), the first cell conducting in arm AP will be cell 1, followed by cell 2, then cell 3, then cell 4, and then back to cell 1 until the state machine transitions to state 1 (S1). In state 1, the first cell conducting in arm AP will be cell 4, then cell 1, etc. The state machine transitions to the next state after a period of $T/2$; since the experiment had an AC period of 200 us, the state machine therefore switches state every 100 us. This is done so that each arm can change its switching order while it is not conducting, preventing unwanted switching behavior. Note that the states alternate which arm has its switching order changed: from S0 to S1, arm AP has its switching order changed while arm AN retains its switching order, then from S1 to S2 arm AP retains its switching order while arm AN has its switching order changed, etc.

Table A-2
Cell-Swapping Algorithm Description

| Order | S0 | S1 | S2 | S3 | S4 | S5 | S6 | S7 |
|---|---|---|---|---|---|---|---|---|
| 1 | AP1 | AP4 | AP4 | AP3 | AP3 | AP2 | AP2 | AP1 |
| 2 | AP2 | AP1 | AP1 | AP4 | AP4 | AP3 | AP3 | AP2 |
| 3 | AP3 | AP2 | AP2 | AP1 | AP1 | AP4 | AP4 | AP3 |
| 4 | AP4 | AP3 | AP3 | AP2 | AP2 | AP1 | AP1 | AP4 |
| 1 | AN4 | AN4 | AN3 | AN3 | AN2 | AN2 | AN1 | AN1 |
| 2 | AN1 | AN1 | AN4 | AN4 | AN3 | AN3 | AN2 | AN2 |
| 3 | AN2 | AN2 | AN1 | AN1 | AN4 | AN4 | AN3 | AN3 |
| 4 | AN3 | AN3 | AN2 | AN2 | AN1 | AN1 | AN4 | AN4 |

## Nomenclature

| | |
|---|---|
| $v_c$* | Arm capacitor voltage reference |
| $d_{dc}$ | DC or common-mode control variable |
| $d_{ac}$ | AC or differential-mode control variable |
| $m$ | Modulation reference, comprised of $d_{dc}$ and $d_{ac}$ |
| $v_{cn}$ | Voltage of nth capacitor |
| $v_{xy}$, $i_{xy}$ | Arm voltage and current: $x = a, b; y = p, n$ |
| $i_{dc}$, $i_{ac}$ | DC and AC currents |
| $R_{arm}$, $L_{arm}$ | Parasitic arm resistance and inductance |
| $N$ | Number of cells per arm |
| $R_{dc}$. $L_{dc}$ | DC cable resistance and inductance |
| $R_{ac}$, $L_{ac}$ | AC load resistance and inductance |
| $R_{dc}'$, $L_{dc}'$ | Equivalent DC resistance and inductance |
| $R_{ac}'$, $L_{ac}'$ | Equivalent AC resistance and inductance |
| $C$, $C_{eq}$ | Cell capacitance & equivalent arm capacitance |
| $v_{\Sigma}$ | Half sum of phase-leg arm capacitor voltages |
| $v_{\Delta}$ | Half diff. of phase-leg arm capacitor voltages |
| $d_{xy}$ | Arm modulation reference: $x = a, b; y = p, n$ |
| $v_{cxy}$ | Arm capacitor voltage: $x = a, b; y = p, n$ |
| $x(k)$ | Sampled variable at time $k$ |
| $x(k + 1)$ | Future value of a variable at time $k + 1$ |

For variables $d_{xy}$ and $v_{cxy}$, the '$x$' is a placeholder for the label of the phase leg ('$a$' or '$b$'), and '$y$' is a placeholder for whether it is for the positive or negative arm ('$p$' or '$n$').

## IX. Biographies

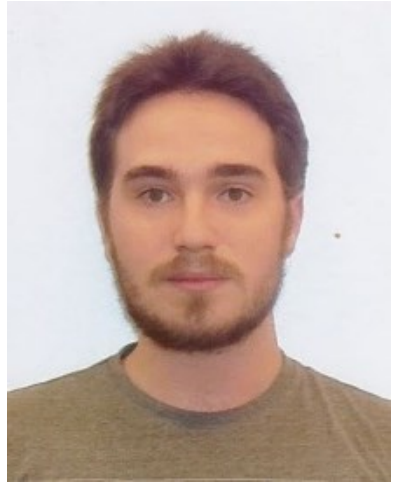

**Sandro Martin** (S'17-M'20) received the B.S. and Ph.D. degrees in Electrical Engineering from the Florida State University, Tallahassee, Florida, USA. Since May 2016, he has been with the Center for Advanced Power Systems at the Florida State University. From May 2018 to September 2018, he was an electrical engineering intern at Mainstream Engineering Co. His special fields of interest include modeling and control, modular multilevel converters, and DC solid-state transformers.

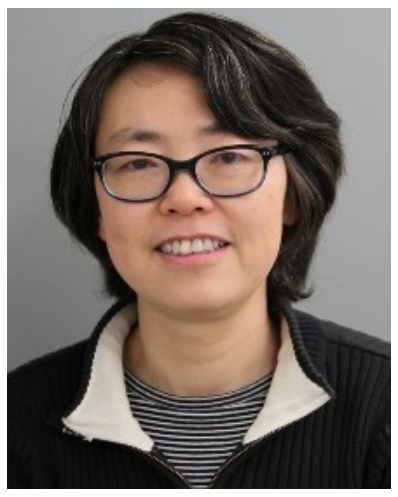

**Hui Li** (S'97–M'00–SM'01) received the B.S. and M.S. degrees in electrical engineering from the Huazhong University of Science and Technology, China, in 1992 and 1995, respectively, and the Ph.D. degree in electrical engineering from the University of Tennessee, Knoxville, in 2000. She is currently a Professor with the Electrical and Computer Engineering Department, Florida A&M University—Florida State University College of Engineering, Tallahassee, FL, USA. Her research interests include PV converters, energy storage applications, and smart grid.

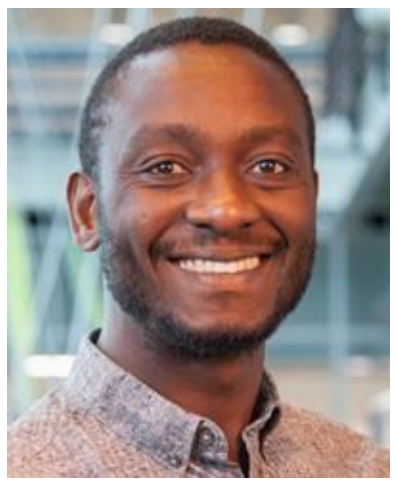

**Olugbenga Moses Anubi** (M'15) received the Ph.D. degree in mechanical engineering from the University of Florida, Gainesville, FL, USA, in 2013. He was a Lead Control Systems Engineer with GE Global Research, Niskayuna, NY, USA. He is currently an Assistant Professor of Electrical and Computer Engineering with the Florida A&M University-Florida State University (FAMU-FSU) College of Engineering, Florida State University, Tallahassee, FL, with affiliations with the Center for Advanced Power Systems and the Center for Intelligent Systems, Controls and Robotics. His research interests include robust, resilient, and adaptive control systems, cyber-physical systems control, real-time optimization, robotics, and vehicle dynamics and control.